\documentclass[11pt]{article}
\usepackage{bridges}
\usepackage{amsfonts,amssymb,amsthm,eucal,amsmath}
\usepackage{graphicx}
\usepackage{subfig}
\usepackage{caption}

\usepackage{wrapfig}
\usepackage{pinlabel, color}
\usepackage{multirow}

\usepackage{hyperref}

\newcommand{\R}{\mathbb{R}}

\newcommand{\Ha}{\mathbb{H}}

\title{The Quaternion Group as a Symmetry Group}
\author{
\begin{tabular}{cc}
Vi Hart            & Henry Segerman \\
Communications Design Group    & Department of Mathematics \\
SAP Labs       & Oklahoma State University \\
San Francisco, CA, USA        & Stillwater, OK, USA \\
vi@vihart.com & henry@segerman.org\\
vihart.com  & segerman.org\\
\end{tabular}}

\date{}

\begin{document}
\maketitle

\begin{abstract}
We briefly review the distinction between abstract groups and symmetry groups of objects, and discuss the question of which groups have appeared as the symmetry groups of physical objects. To our knowledge, the quaternion group (a beautiful group with eight elements) has not appeared in this fashion. We describe the quaternion group, both formally and intuitively, and give our strategy for representing the quaternion group as the symmetry group of a physical sculpture.
 \end{abstract}

\section{Introduction}
\label{Sec:Intro}

A \emph{symmetry} of an object is a geometric transformation which leaves the object unchanged. So, for example, an object with 3-fold rotational symmetry has three symmetries: rotation by $120^\circ$, rotation by $240^\circ$, and the trivial symmetry, where we do nothing. The symmetries of an object naturally form a group under composition. 
Care must be taken to differentiate between the symmetry group of an object, consisting of geometric transformations that leave the object unchanged, and the abstract group, which only contains information about how the elements of the group interact with each other under composition.

\begin{figure}[htbp]
\centering 

\subfloat[An object with symmetry group isomorphic to $D_4$.]
{
\includegraphics[width=0.3\textwidth]{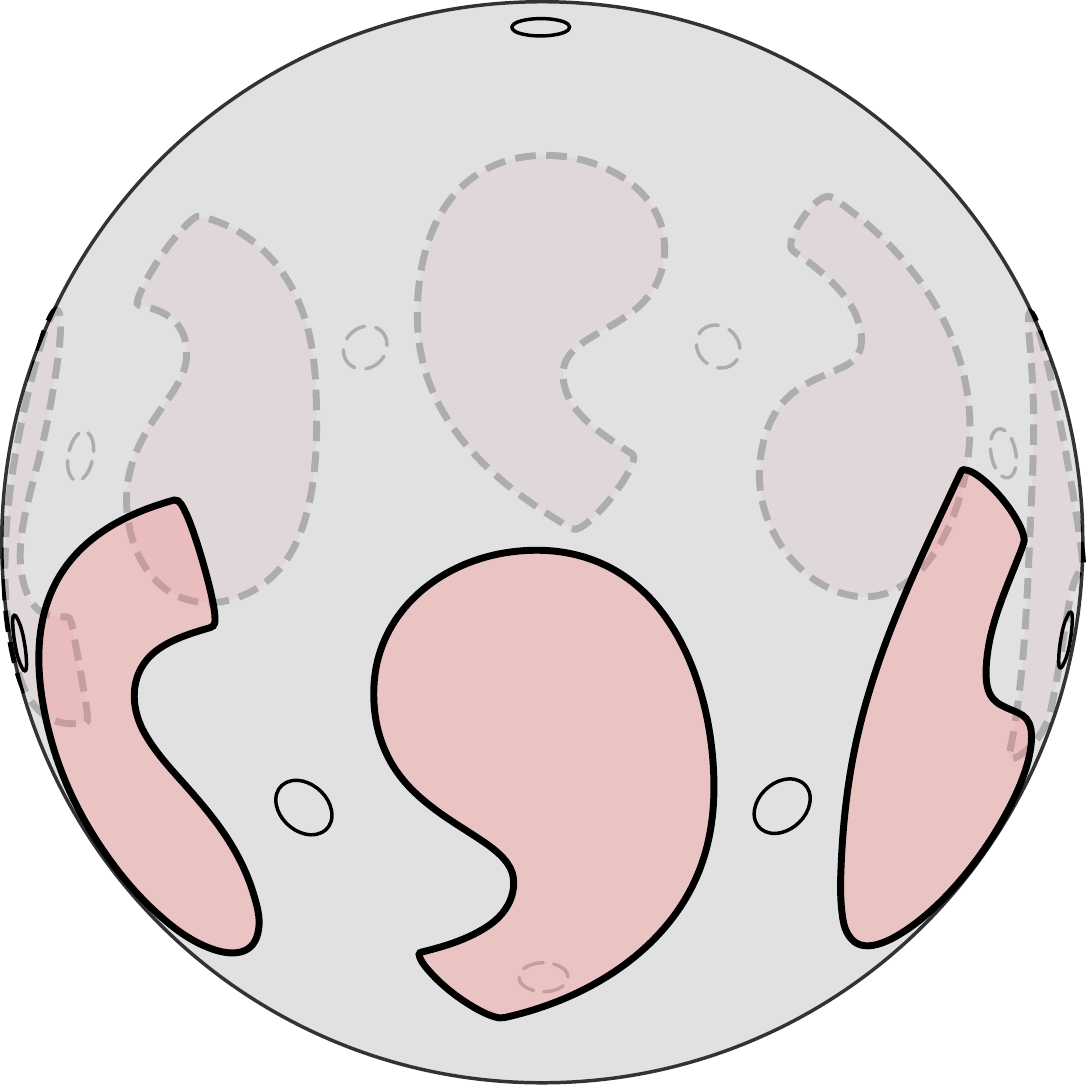}
\label{Fig:224}
} 
\quad
\subfloat[Another object with symmetry group isomorphic to $D_4$.]
{
\includegraphics[width=0.3\textwidth]{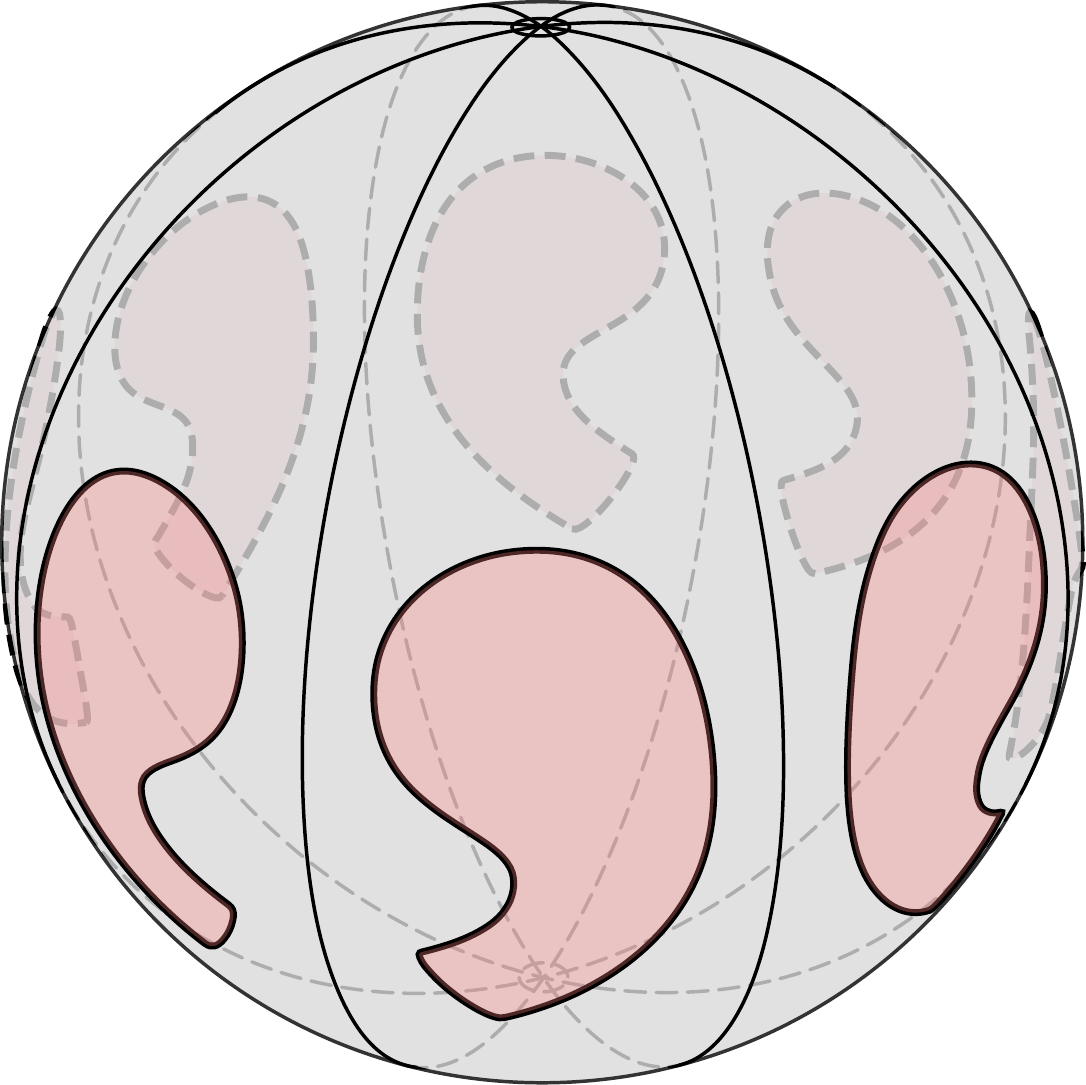}
\label{Fig:*44}
}
\quad
\subfloat[A Cayley graph for $D_4$. The edges with arrows correspond to rotations, the other edges correspond to reflections.]
{
\includegraphics[width=0.3\textwidth]{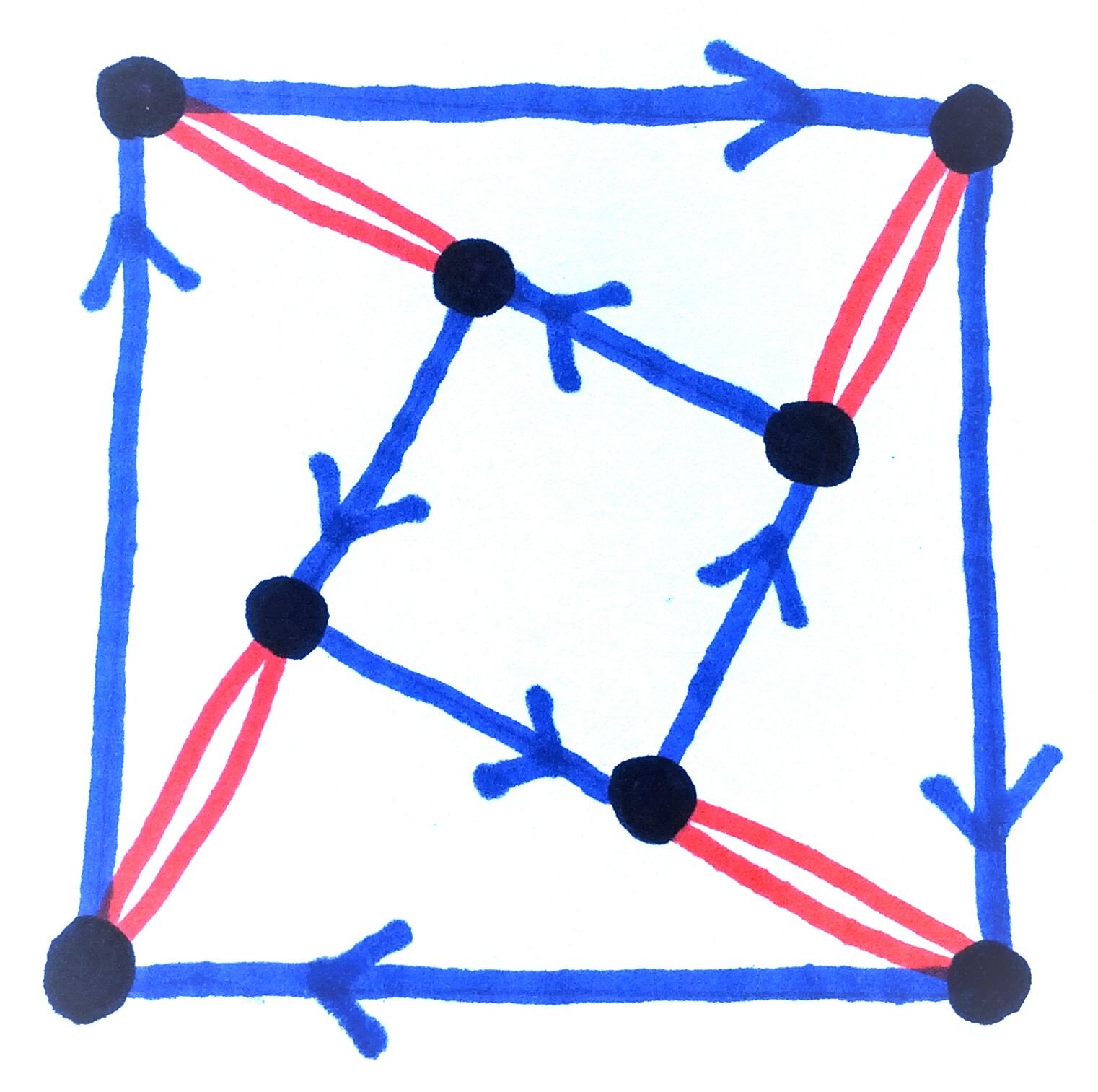}
\label{Fig:D4_cayley}
}

\caption{Symmetric designs on the sphere.}
\label{Fig:spherical_symm_exs}
\end{figure}

As an example, consider the objects pictured in Figure \ref{Fig:spherical_symm_exs}. The object shown in Figure \ref{Fig:224} has no planes of mirror symmetry, but has ten \emph{gyration points} (points of rotational symmetry, marked with a small circle). 
The symmetry group of this object consists of rotations about these gyration points by multiples of either $90^\circ$ or $180^\circ$ (depending on the kind of gyration point).
The object shown in Figure \ref{Fig:*44} has four planes of mirror symmetry, marked by their intersection with the sphere. The symmetry group of this object consists of reflections in these planes, together with rotations about the gyration points at the top and bottom of the sphere.

These objects have different sets of symmetries, although as abstract groups they are both isomorphic to the dihedral group with eight elements, $D_4 = \langle a,b \mid a^{4} = b^2 = 1, b a b = a^{-1}\rangle$. 
Figure \ref{Fig:D4_cayley} shows a Cayley graph which abstractly represents the group $D_4$, although the figure itself does not have $D_4$ symmetry.\footnote{Cayley graphs are purely combinatorial objects with no intrinsic embedding in space associated with them, although they are often drawn symmetrically embedded in space (in this case with order 4 rotational symmetry). Under our definition then, the Cayley graph of a group does not count as a physical object with the symmetry of that group, unless it is embedded in a way that actually has that symmetry.}

As artists, we are particularly interested in 
the symmetries of real world physical objects.
Three natural questions arise:
\begin{enumerate}
\item Which groups can be represented as the group of symmetries of some real-world physical object? 
\item Which groups \emph{have actually} been represented as the group of symmetries of some real-world physical object?
\item Are there any glaring gaps -- small, beautiful groups that should have a physical representation in a symmetric object but up until now have not?
\end{enumerate}

In one sense, the first question is answered in detail in \emph{The Symmetries of Things} \cite{symmetries_of_things} by Conway, Burgiel and Goodman--Strauss. For example, the spherical symmetries of objects in $\R^3$ with symmetries consisting of rotations and reflections whose axes and reflection planes all pass through a single point are entirely classified. 
In another sense things are not so clear -- we can draw pictures of the wallpaper symmetry groups, but real, physical wallpaper never extends infinitely far in the plane. Therefore the design is not actually unchanged by a translation, or even approximately unchanged. However, the viewer can infer how the pattern should continue, and so in this sense the wallpaper groups can be and have been represented by real world objects. Similarly, no sculpture (that we know of) is exactly symmetric at a subatomic level, but the viewer infers the intended symmetry. Using models such as the Poincar\'{e} disk model of hyperbolic space we can make physical representations of hyperbolic symmetry groups, for example in Escher's ``Circle Limit'' series. Likewise, artists have used shadows and projections to represent higher dimensional symmetry groups as physical objects. 

To our knowledge, the second question has not been investigated in this generality.\footnote{A related question that \emph{has} been investigated in detail is the issue of which of the 17 wallpaper groups appear in the Alhambra. This is a topic of some debate (see for example \cite{grunbaum_alhambra_symm_groups}), with questions over whether or not different colours can be counted as the same in order to show a symmetry group, how large a design needs to be to represent a group and so on. However, whether or not these symmetry groups appear in the Alhambra, they certainly have appeared in modern works.}
There are of course infinitely many groups, but only finitely many physical objects that have ever been made. So for example, there is some smallest finite N for which there has never been a physical object with N-fold rotational symmetry. 
Certainly all of the polyhedral symmetry groups have been represented in sculpture since antiquity, if not represented in nature beforehand. As we move to representations of more esoteric groups (for example requiring a model of hyperbolic space), the number of physical representations surely goes down. 

As for the third question, there does seem to be such a gap for the \emph{quaternion group} $Q_8$.\footnote{Note that there are many drawings of the Cayley graph of $Q_8$ with various embeddings, as well as art representing its multiplication table. See for example Gwen Fisher's Quaternion Quilt \cite{fisher}.}  
In Section \ref{Sec:Quat group} we give an abstract definition of the quaternion group. 
In Section \ref{Sec:Visualise quat group} we give an intuitive way to visualise the action of $Q_8$ on a hypercube, by building the boundary of the hypercube from eight decorated cubes. 
In Section \ref{quat group in sculpture} we  
give our strategy for producing a sculpture whose symmetry group is $Q_8$, using stereographic projection from the unit sphere in 4-dimensional space. To our knowledge, this is the first sculpture that specifically exhibits $Q_8$ as a symmetry group.
In Section \ref{future directions} we give some directions for future work.

\section{The Quaternion group}
\label{Sec:Quat group}

The quaternion group $Q_8$ is a beautiful group of order eight, not to be confused with the real quaternions or the unit quaternions.
The real quaternions $\Ha$ are a number system that extends the real numbers, similarly to the way in which the complex numbers also extend the real numbers.\footnote{The symbol $\Ha$ is chosen in honour of the discoverer of the real quaternions, William Rowan Hamilton.} 
$$\Ha = \{w +xi+yj+zk \mid w,x,y,z\in\R, i^2=j^2=k^2=ijk=-1 \}$$
The quaternion group consists of eight elements of $\Ha$ under multiplication.
$$Q_8=\{\pm1, \pm i, \pm j, \pm k\}$$

\begin{wraptable}[17]{r}{0.46\textwidth}   
\vspace{-15pt}
\centering
\[
\begin{array}{ r | rrrrrrrr }

  & 1 & -1 & i & -i & j & -j & k & -k \\
  \hline
1 & 1 & -1 & i & -i & j & -j & k & -k \\
-1 & -1 & 1 & -i & i & -j & j & -k & k \\
i & i & -i & -1 & 1 & k & -k & -j & j\\
-i & -i & i & 1 & -1 & -k & k & j & -j\\
j & j & -j & -k & k & -1 & 1 & i & -i\\
-j & -j & j & k & -k & 1 & -1 & -i & i\\
k & k & -k & j & -j & -i & i & -1 & 1\\
-k & -k & k & -j & j & i & -i & 1 & -1

\end{array}
\]

\caption{The multiplication table (i.e., Cayley group table) for the quaternion group. The convention for the order of multiplication used here is that the row label is first, followed by the column label.}
\label{mult_table}
\end{wraptable}

We give the multiplication table for $Q_8$ in Table \ref{mult_table}. The quaternion group has one element of order one, one element of order two, and  six elements of order four. 
We can view $\Ha$ as $\R^4$ with an added multiplication structure. Multiplication by a quaternion then acts as a transformation on $\R^4$. The unit quaternions are quaternions $q=w+xi+yj+zk$ with magnitude $|q|=\sqrt{w^2+x^2+y^2+z^2}=1.$ These act by 4-dimensional rotation on $\R^4$, and so $Q_8$ acts by rotation on $\R^4$ (and also acts by rotation on $S^3$, the sphere of radius 1 in $\R^4$). In fact, $Q_8$ also acts nicely on the eight cubical cells of the hypercube, permuting them. To see this, note that the eight cells are given by $\{(1,x,y,z)\mid x,y,z\in\R\}$, $\{(-1,x,y,z)\mid x,y,z\in\R\},$ and similarly but with the $1$ or $-1$ in the three other possible coordinate positions. If, for example, we act on the point $(1,x,y,z) = 1+ix+jy+kz$ on the right by $i$, we get $(1+ix+jy+kz)i = i-x-ky+jz=-x+i+zj-yk=(-x,1,z,-y)$, which is on one of the other cells.

\section{Visualizing the quaternion group}
\label{Sec:Visualise quat group}

\begin{wrapfigure}[14]{l}{0.35\textwidth}
\vspace{-10pt}
\centering 

\includegraphics[width=0.35\textwidth]{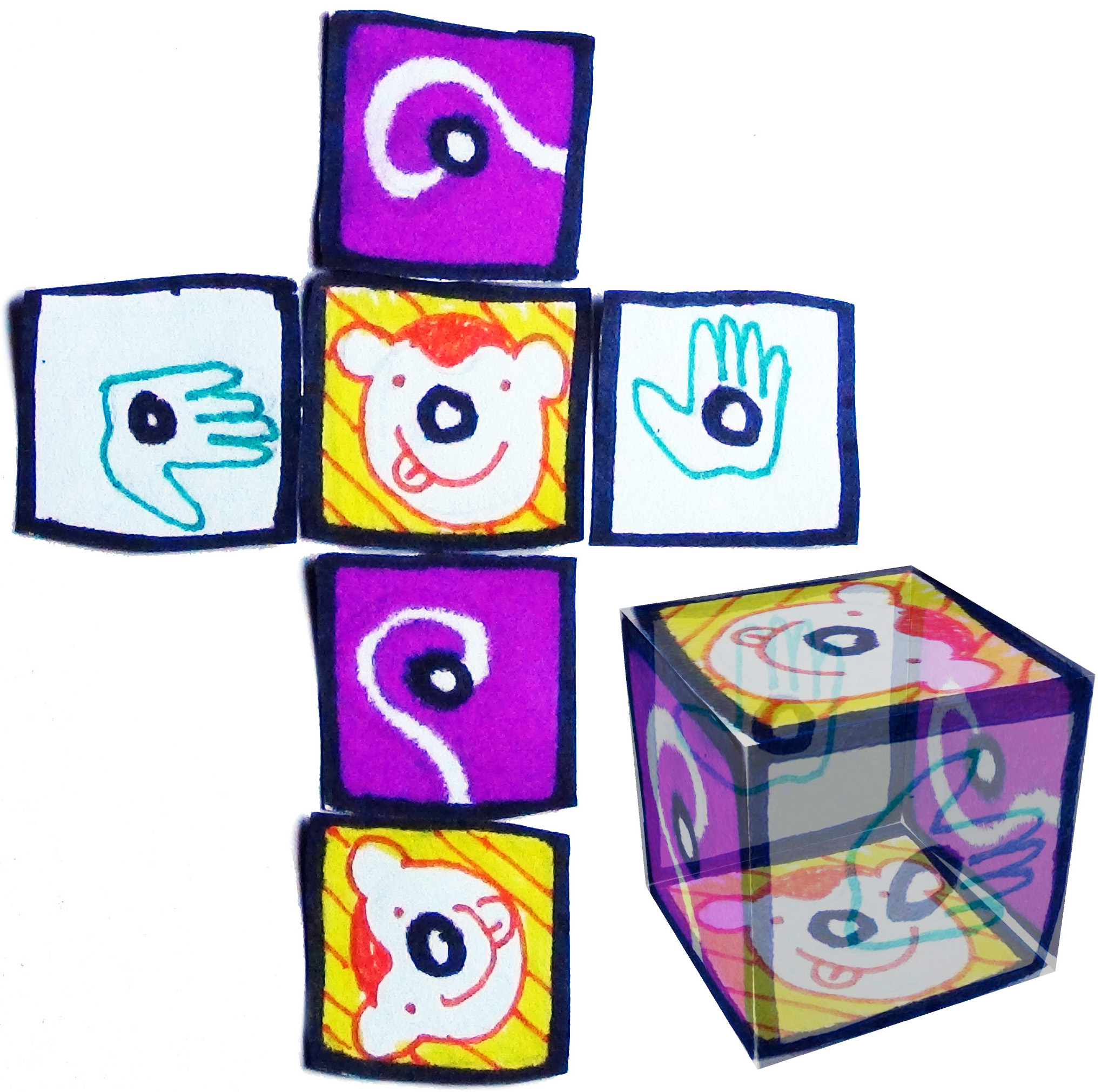}
\caption{A monkey block.} 
\label{monkeyblock_net}
\end{wrapfigure}

In order to visualize the symmetries of the quaternion group, we begin by assembling a hypercube out of eight \emph{monkey blocks}.\footnote{The reason for our use of a monkey will become clear in Section \ref{quat group in sculpture}.} See Figure \ref{monkeyblock_net}. Each monkey block is a cube which, due to the patterns on the faces, has no symmetry itself. One face of the cube has a monkey's face on it with her tongue sticking out to the right, while the opposite face on the cube has a monkey face with her tongue sticking out to the left. In addition, opposite faces are rotated a quarter turn from each other. 
A right monkey paw is opposite a quarter-turned left monkey paw, and a question-mark tail is opposite a quarter-turned un-question-mark tail.

\begin{figure}[h]
\centering 

\includegraphics[width=0.8\textwidth]{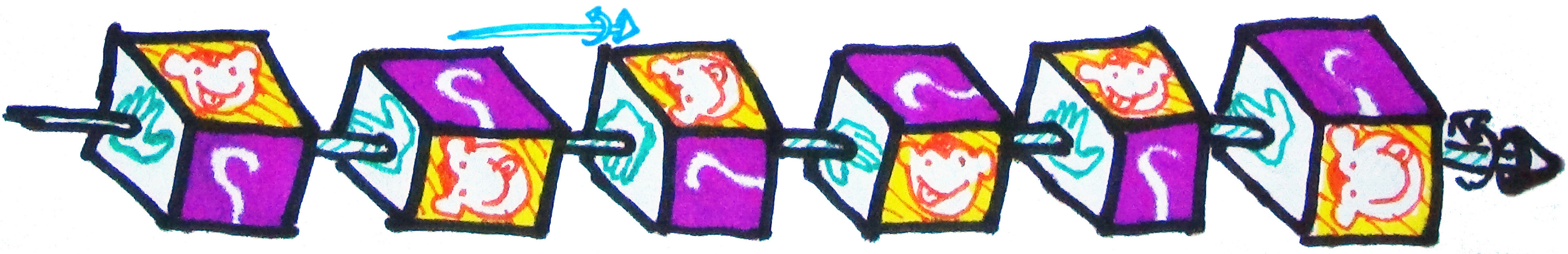}

\caption{An infinite line of monkey blocks.} 
\label{line}
\end{figure}

Given two of these monkey blocks, one can match any face on one block to its mirrored version on the other. Given an infinite number, one can match up all the right paws with all the left paws to get an infinite line of monkey blocks. See Figure \ref{line}. This infinite line has translational symmetry: the whole line can move forward or backwards any multiple of four blocks. It also has a more unusual symmetry: it can move forward one and rotate $90^\circ$ to the left, giving a screw motion. 

\begin{wrapfigure}[13]{r}{0.3\textwidth}
\vspace{-10pt}
\centering 
\includegraphics[width=0.3\textwidth]{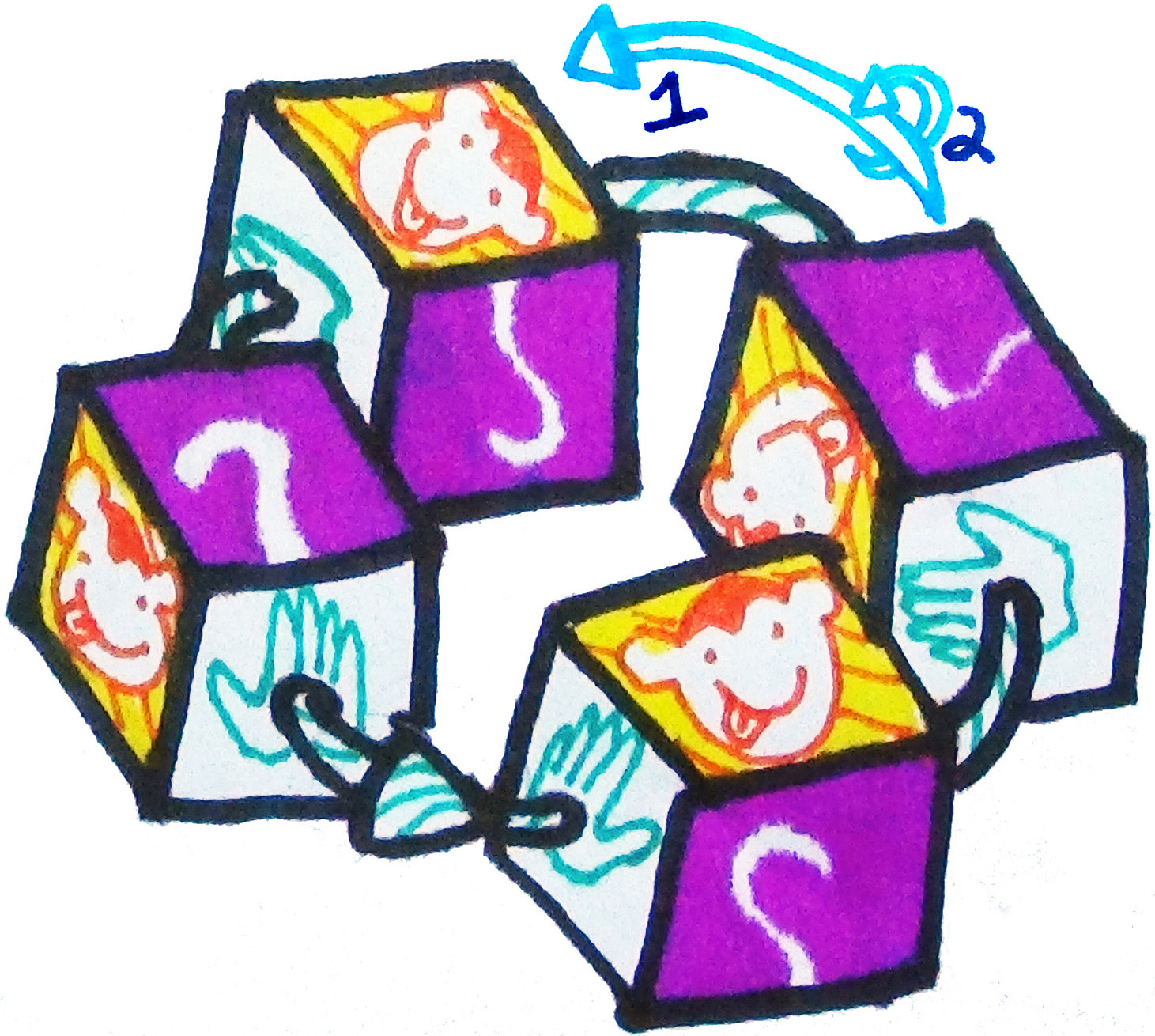} 
\caption{A ring of four monkey blocks.} 
\label{ring}
\end{wrapfigure}

The next obvious thing to try when matching up cubes might be tiling monkey blocks infinitely in two or three dimensions, but there is no way to match four blocks around a single edge.
However, if we allow ourselves to attach blocks with a $90^\circ$ bend into the fourth dimension, we can take four blocks from our infinite line and attach them in a loop. See Figure \ref{ring}. 
In this four-dimensional arrangement (which is not really shown in Figure \ref{ring}!) we have two perpendicular planes of rotation. Rotation in the first plane corresponds to moving the blocks around the circle passing through the paw sides of the monkey blocks in the direction of the arrow (labelled with a ``1''). Rotation in the second plane (which is harder to see in the figure
) corresponds to keeping each block in place but spinning it on the circle (the arrow labelled ``2'').

The symmetry of the ring of blocks consists of simultaneous\footnote{Two perpendicular rotations in four dimensions may be performed in either order to obtain the same result.} $90^\circ$ rotation in the two perpendicular planes, something not possible in fewer than four dimensions. This makes the symmetry an \emph{isoclinic rotation}: two equal and perpendicular rotations. Isoclinic rotations are in fact the only symmetry of this ring of monkey blocks --
it has none of the symmetries found in three dimensions.

Those familiar with four dimensional objects may notice that this ring of four cubes is half of the boundary of a hypercube. It might seem unlikely that by continuing to match the same monkey blocks face to face one could complete the hypercube without violating the matching rules. However, 
not only does this work, but 
the set of eight matched monkey blocks have symmetry around all three possible sets of isoclinic rotations: the $wx$-plane with the $yz$-plane, the $wy$-plane with the $xz$-plane, and the $wz$-plane with the $xy$-plane.
There are no other symmetries.

\begin{figure}[htbp]
\centering 
\hspace{-.7cm}
\subfloat[Eight monkey blocks arranged into a hypercube net.]
{
\includegraphics[width=0.48\textwidth]{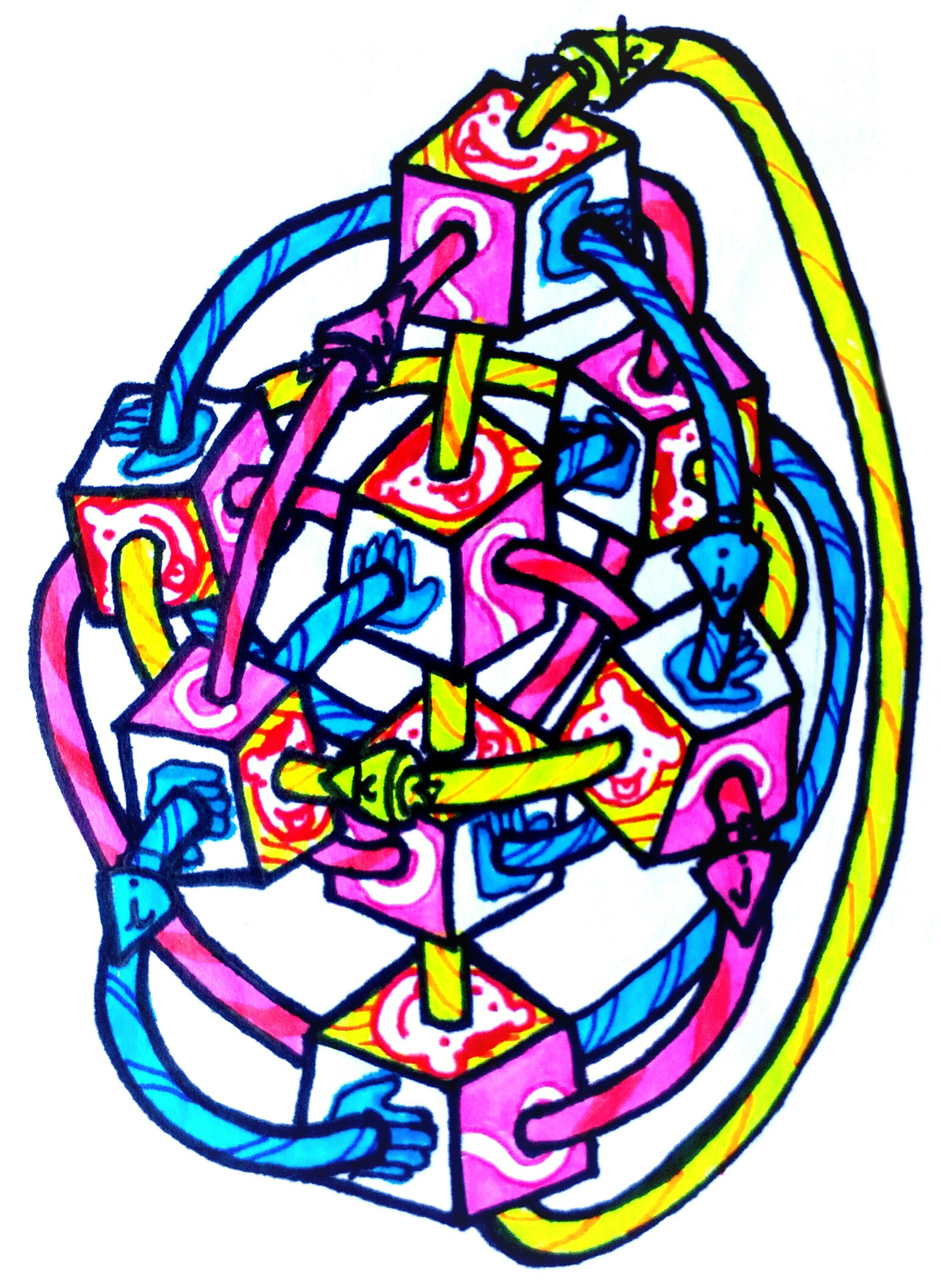}
\label{monkeyblock_hypercube_net}
}
\hspace{.7cm}
\subfloat[The hypercube net labelled with the elements of $Q_8$. This can be interpreted as the Cayley graph for $Q_8$ with generating set $\{i,j,k\}$.]
{
\includegraphics[width=0.43\textwidth]{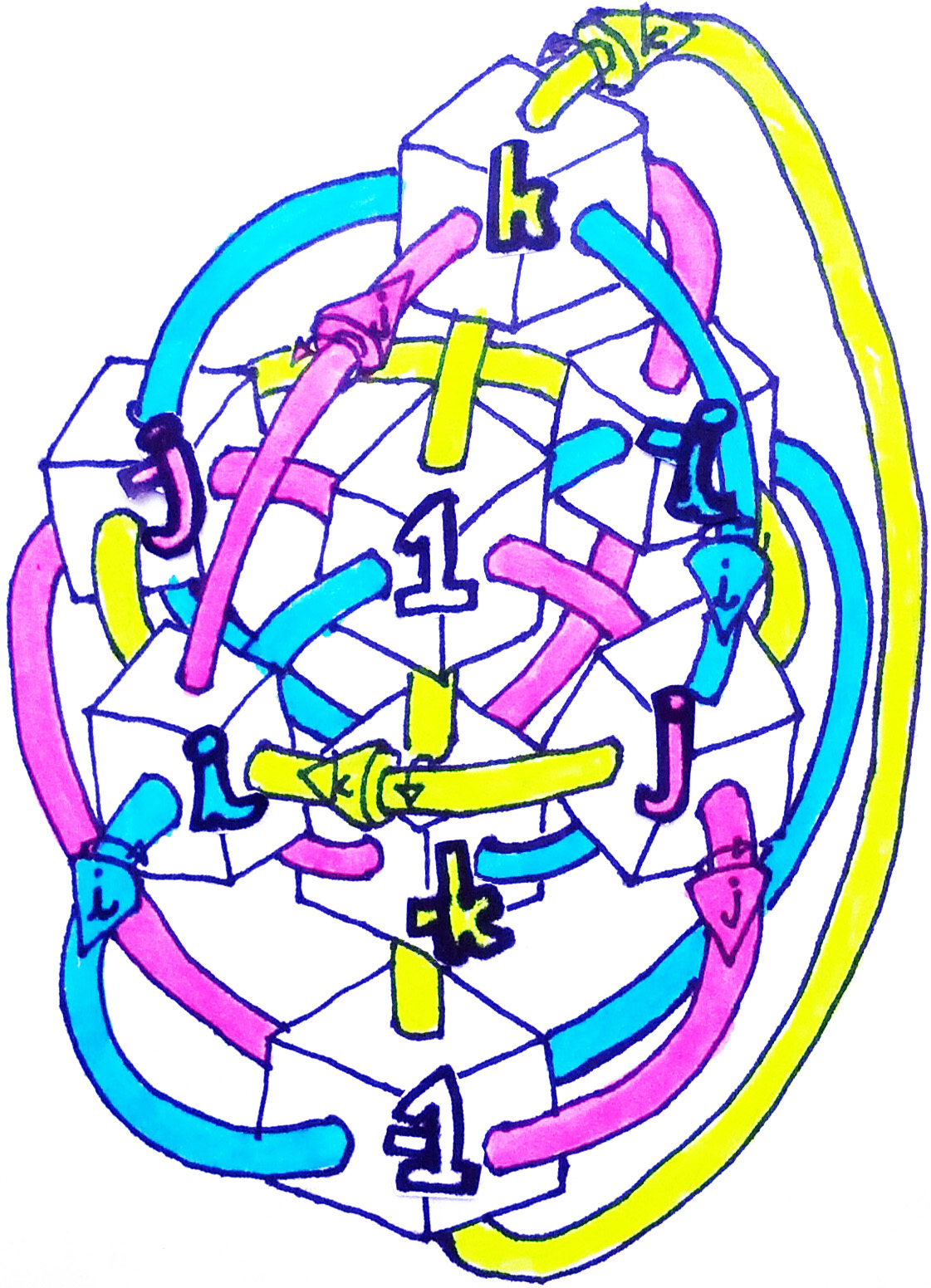}
\label{elements}
} 
\caption{}
\label{monkey block hypercubes}
\end{figure}

In Figure \ref{monkeyblock_hypercube_net} the monkey blocks are shown as the net of a hypercube, with the three isoclinic rotations shown as pairs of axes. The monkey tails rotate around in two perpendicular loops of four blocks each, which we label as the $i$ axes. The $j$ axes are the two loops of four blocks made by connecting matching hands, and the $k$ axes go through the face faces. The $k$ axis is perhaps the easiest to see as the four blocks in one loop are arranged vertically, spinning up one step and a quarter turn to the left each time. At the same time, in the perpendicular plane the other four blocks move likewise: a quarter turn to the left around the vertical $k$ axis while rotating a quarter turn to the left around the horizontal $k$ axis.

Our labelling of the axes suggests that these isoclinic rotations can be thought of as multiplication by the quaternions $i$, $j$, and $k$, and of course this is true. Starting from an arbitrarily chosen block which we label 1, we can label each of the eight cells of the hypercube by the element of $Q_8$ we need to multiply (on the right) to get there. See Figure \ref{elements}. 

Starting with any labelled cube, if we follow the arrows and write down the sequence of labels, then we finish at the cube labelled with the element of $Q_8$ which is the product of the starting cube label with the arrow labels, read in order. Following an arrow backwards corresponds to multiplying by the negative of the label on the arrow. Using these rules, we can see that rotating twice along any of the three axes (an isoclinic rotation of $180^\circ$) will take us to the opposite cube, as will rotating by $90^\circ$ once each along the three axes in order, geometrically showing that $i^2 = j^2 = k^2 = ijk = -1$. 

\section{The quaternion group in sculpture}
\label{quat group in sculpture}

While the quaternion group is well known as an abstract group, as are its actions on objects such as the hypercube and the 3-sphere, to our knowledge no physical object has been made which has as its symmetry group the quaternion group.

It is a reasonably simple exercise to check that $Q_8$ is not a subgroup of the isometries of 3-dimensional Euclidean space. However, as noted in Section \ref{Sec:Quat group}, $Q_8$ is very naturally a subgroup of the symmetries of the 3-dimensional sphere, $S^3$, the unit sphere in $\R^4$. Using techniques we described previously in \cite{sculptures_in_S3}, we can create a sculpture that has $Q_8$ as its symmetry group. As noted in Section \ref{Sec:Quat group}, the eight elements of $Q_8$ act on the eight cubical cells of the hypercube. We radially project the boundary of the hypercube onto $S^3$, then use stereographic projection to move the eight cubes into $\R^3$, where we can build a physical sculpture using 3D printing. We must choose a design to put in each of the eight cubes. The design must itself have no symmetry, or else the symmetry group would be larger than $Q_8$. In order to be produced as a physical sculpture, the design in each cubical cell of the hypercube also needs to physically connect with its copies in the six neighbouring cubical cells, through the six square faces of the cube. The design must connect onto copies of itself following the same matching rules as the monkey blocks of Section \ref{Sec:Visualise quat group}. 

\begin{wrapfigure}[26]{r}{0.63\textwidth}
\vspace{-5pt}
\centering 

\includegraphics[width=0.63\textwidth]{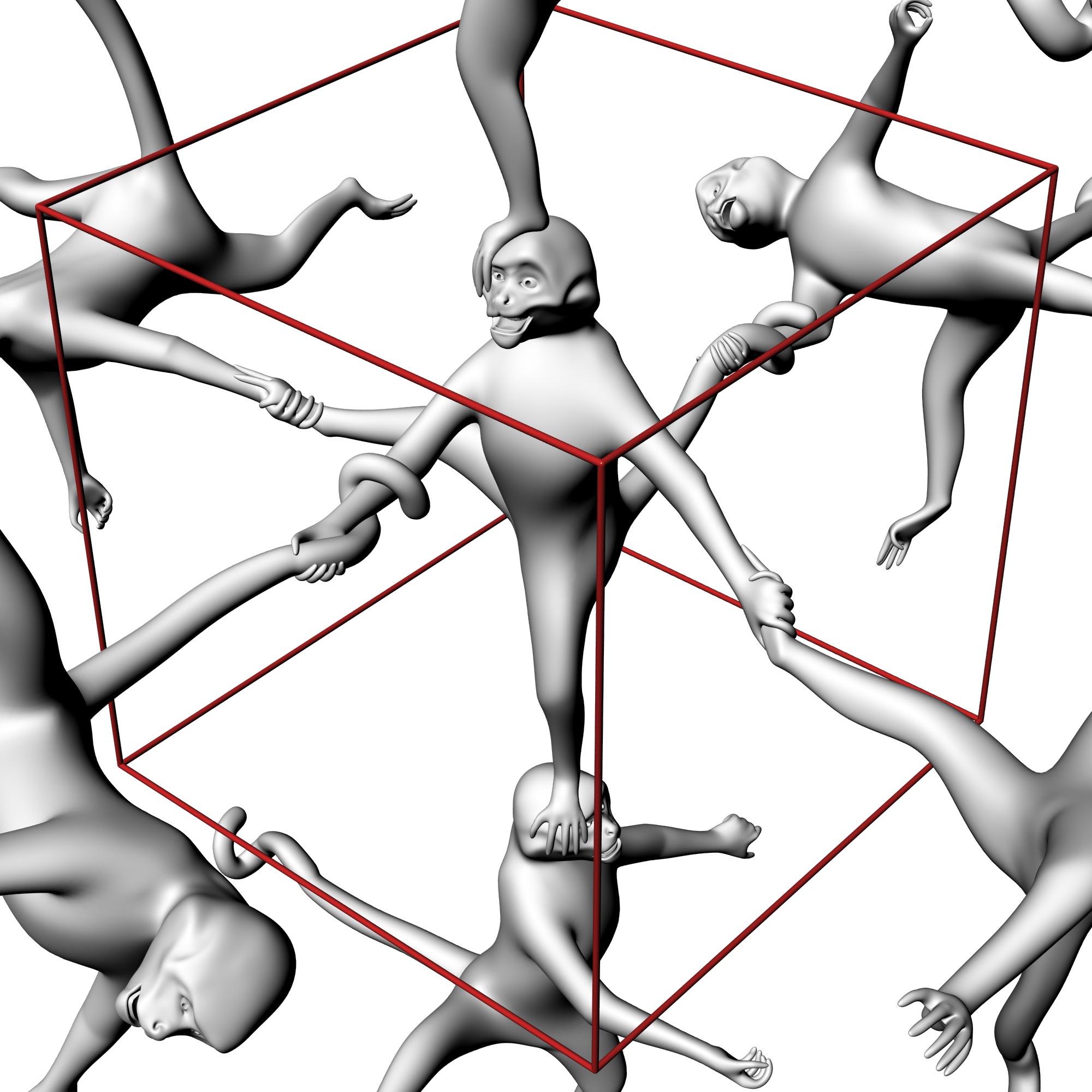}

\caption{The monkey in a Euclidean cube, prior to mapping into $S^3$, with six neighbours.} 
\label{euclidean_monkeys_in_cube}
\end{wrapfigure}

The obvious (perhaps even canonical) choice for such a design is, of course, a monkey. See Figure \ref{euclidean_monkeys_in_cube}. With appropriate posing of the monkey, its bilateral symmetry can be broken, and including the head and tail it has six limbs, one for each face of the cube. The monkey's left foot stands on the head of a neighbour, the left hand grabs a neighbour's right foot, and the right hand grabs a neighbour's tail. By symmetry, everything that goes around comes around -- so the other three neighbours of this monkey are standing on this monkey's head, grabbing its right foot, and grabbing its tail.

The monkey was designed in a Euclidean cube. It was then run through eight different transformations in order to move eight copies of it to the appropriate positions in $S^3$ and then back to $\R^3$ by stereographic projection. The first step of all of these transformations is to project the Euclidean cube into a curved cube in $S^3$. This is done in exactly the same way that a cubical cell of the hypercube is radially projected onto the hypersphere $S^3$. To be precise, we think of a point $(x,y,z)$ in the Euclidean cube $[-1,1]^3$ as actually being the point $(1,x,y,z)$ on one of the cells of the Euclidean hypercube $[-1,1]^4$, and map it to $S^3$ by

$$  (x,y,z)\mapsto \frac{(1,x,y,z)}{\sqrt{1+x^2+y^2+z^2}}.$$

Now that the design is on $S^3$, we (right) multiply it by $1,i,j,k,-1,-i,-j$ and $-k$ respectively for the eight transformations, and stereographically project each back to $\R^3$. There is a choice of where to project from --  we put the north pole at a vertex of the hypercube, so that in the projection the copies of the monkey are as far from infinity as possible. This makes the resulting features of the eight monkeys as near in size to each other as possible. Very small features may be too fragile to 3D print -- to avoid this we can scale the entire sculpture up, but only so far as the largest features fit within the printer and our budget.

\begin{figure}[hp]
\centering 

\subfloat[]
{
\includegraphics[width=.72\textwidth]{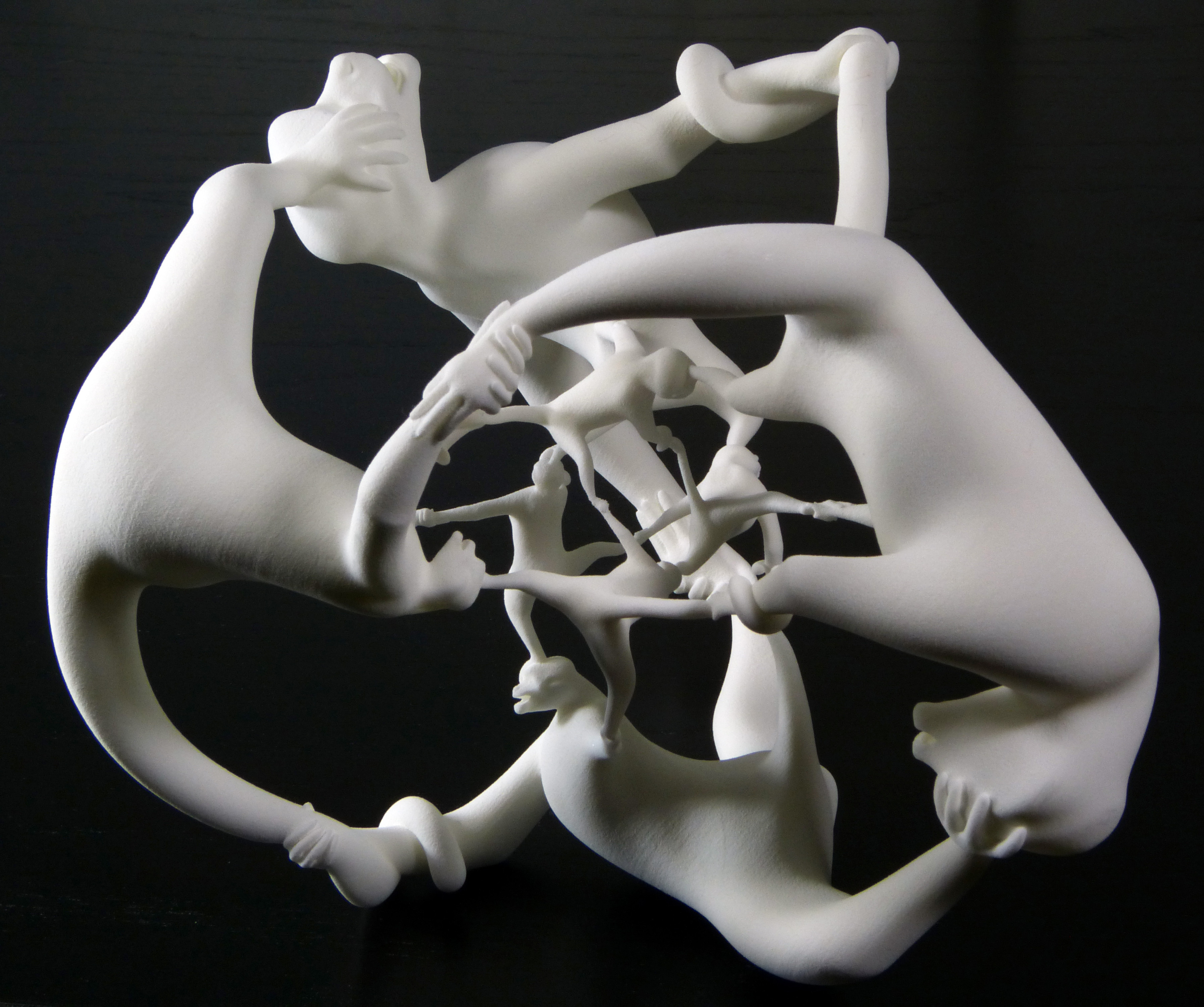}
\label{monkeys_1}
}

\subfloat[]
{
\includegraphics[width=.72\textwidth]{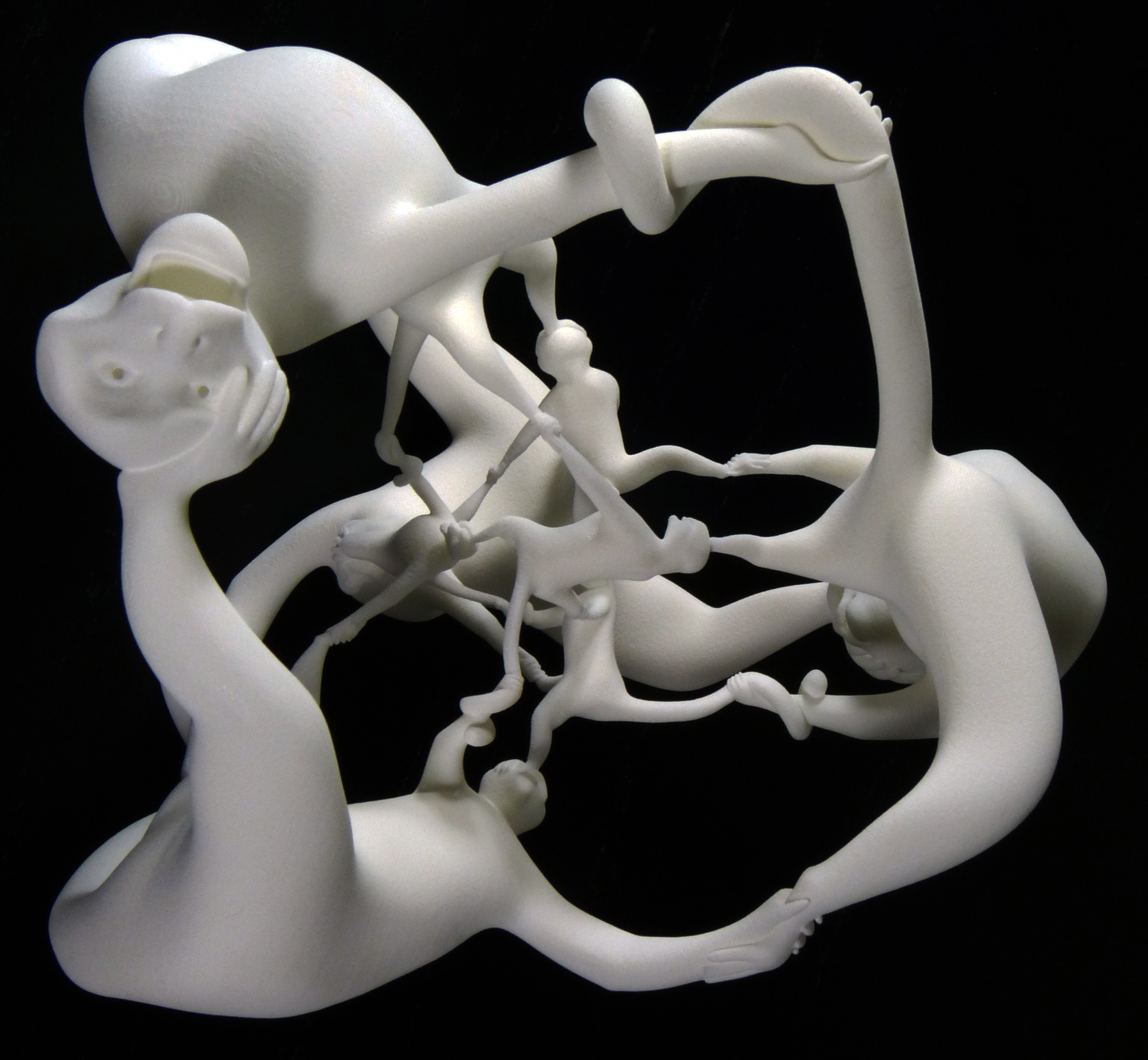}
\label{monkeys_3}
} 
\caption{Two views of the sculpture ``More fun than a hypercube of monkeys''. 
}
\label{fig_monkeys}
\end{figure}

The resulting sculpture is shown in Figure \ref{fig_monkeys}.  Note that the sculpture has no ``ordinary'' symmetries at all: every monkey is different if we only consider isometries of 3-dimensional space. However, under the appropriate isometries of the 3-sphere (as seen through the lens of stereographic projection) they are all identical. From this vantage point the three pairs of axes of rotation have equal billing: each circle consists of two larger, outer monkeys and two smaller, inner monkeys. The three pairs of axes go through the head-foot, hand-foot and hand-tail connections. 

\begin{wrapfigure}[10]{r}{0.3\textwidth}
\vspace{-10pt}
\centering 

\includegraphics[width=0.3\textwidth]{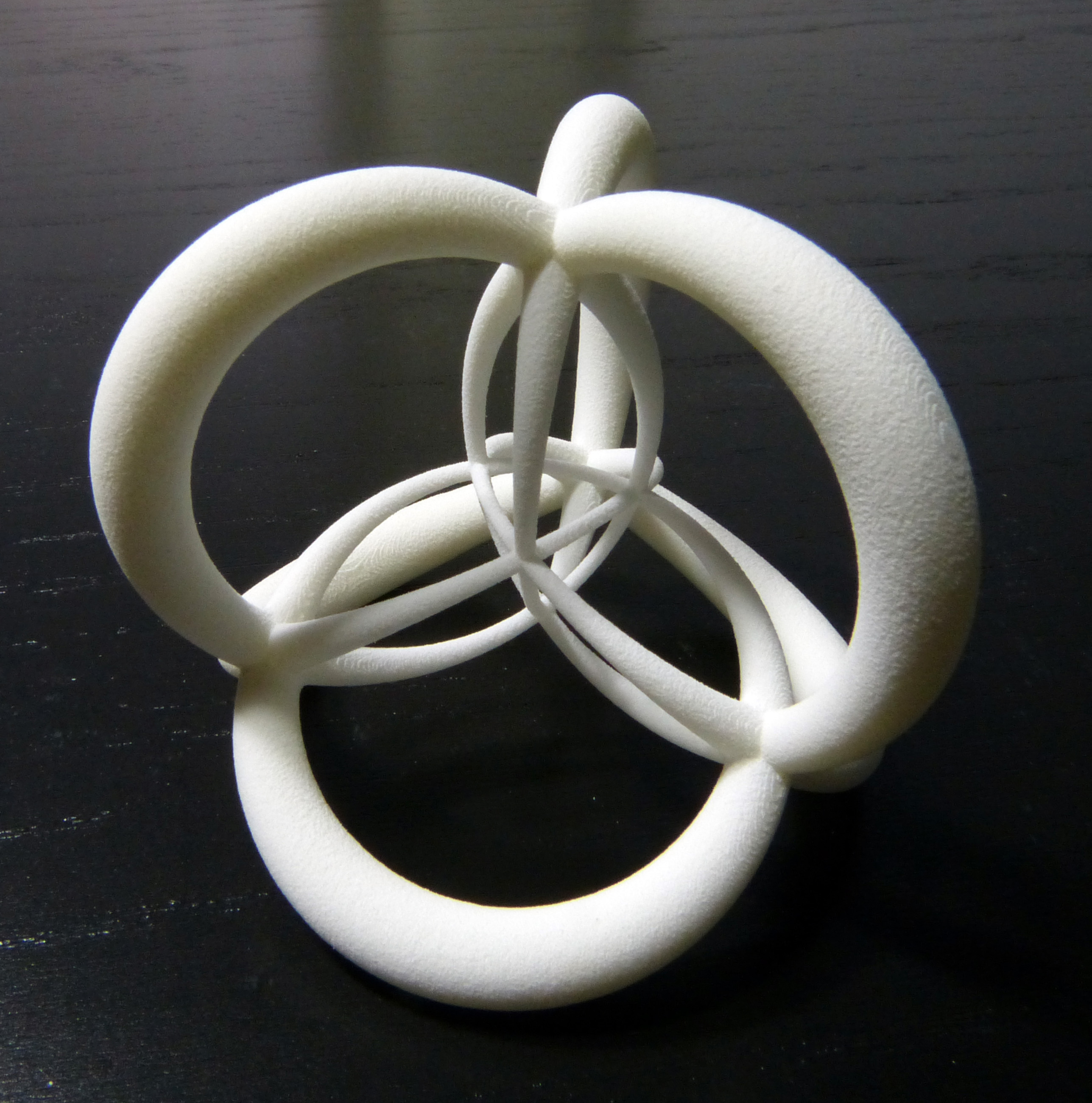}

\caption{The 16-cell.} 
\label{sixteencell} 
\end{wrapfigure}

Each monkey sits inside of a cell of the hypercube and connects to its neighbours through the 2-dimensional faces of the hypercube. Therefore, taken together they form the edges and vertices of the dual polytope to the hypercube, namely the 16-cell. Figure \ref{sixteencell} shows a 3D print of the edges and vertices of the 16-cell, projected to 3-dimensional space in the same way that the monkeys are. This perhaps shows the pairs of axes of rotation more clearly. The symmetry group of this object contains $Q_8$ but is much larger -- the asymmetric monkey design is necessary in order for us to represent \emph{only} the quaternion group.

The object shown in Figure \ref{Fig:*44} is \emph{achiral}, meaning that it is the same as its mirror image. The object in Figure \ref{Fig:224} is \emph{chiral} -- the mirrored object is different.\footnote{To be precise, an object is \emph{achiral} if there is an orientation preserving isometry which moves the mirrored object to coincide with the original. An object is \emph{chiral} if this is only possible with an orientation reversing isometry.} However, the symmetry group of the mirrored object is the same as that of the unmirrored one.\footnote{Or rather, it is a conjugate of the original group of isometries by an orientation preserving isometry.} This is true of almost all objects, if we restrict ourselves to only consider the ordinary symmetries of 3-dimensional objects.\footnote{There are a few metachiral space groups in three dimensions, but no point groups.} However, our monkey sculpture with $Q_8$ symmetry is different. In Conway's classification it is \emph{metachiral} (see \cite{on_quaternions_and_octonions}, p44). Not only is it different from its mirror image, but also all of the left screw symmetries become right screw symmetries in the mirror image, and so the symmetry group itself is also reflected!\footnote{By which we really mean that it is a conjugate of the original group of isometries only by an orientation \emph{reversing} isometry.} In other words, the symmetry group itself has a handedness.

\section{Future directions}
\label{future directions}

It would be an interesting project to list the abstract groups against works of art or illustration, or objects from nature with the given group as its symmetry group.
Are there reasonably simple polytopes (perhaps shaped somewhat less like monkeys) that have symmetry group $Q_8$? Further new sculptures could be produced via stereographic projection, using other finite subgroups of the isometries of $S^3$, as listed in Tables 26.1 through 26.3 of \cite{symmetries_of_things}. 
 Are there other natural strategies for making objects exhibiting previously unrepresented symmetry groups?

\section{Acknowledgements}

We thank Marc ten Bosch for motivating this line of thought with some interesting questions about computer graphics in four dimensions, and Will Segerman for design of the monkey.

\bibliographystyle{hyperplain} 
\bibliography{biblio}

\end{document}